\documentclass[a4paper]{article}
\usepackage{graphicx,epsf,psfrag,tabularx,subfigure,enumerate}
\usepackage[footnotesize,sc]{caption}
\usepackage{amsmath}
\usepackage{multicol}
\usepackage{setspace}
\usepackage[dvips]{color}
\usepackage{soul}
\usepackage{epstopdf}

\renewcommand{\thesubfigure}{\alph{subfigure}}
\makeatletter
  \renewcommand{\@thesubfigure}{(\thesubfigure)\space}
  \renewcommand{\p@subfigure}{\thefigure}
\makeatother

\begin{document}

\begin{center}
{\bf \large On the Gaussian Curvature of Creased Tubes}

K A Seffen \& C R Calladine

\footnotesize{\em Department of Engineering, University of
Cambridge, Trumpington Street, Cambridge, CB2 1PZ.}

\footnotesize{T: +44 1223 764137; F: +44 1223 332662; E:
kas14@cam.ac.uk}

\end{center}

\begin{abstract}
We calculate the Gaussian curvature of a curved, twisted crease in terms of the rate of change of solid angle along its length; we find that this depends on the fold angle across the crease and on the curvature along it but is independent of twist.  We use this result to resolve a paradox concerning the geometry of a creased, twisted-prismatic tube; that there can be no Gaussian curvature overall despite the surface being doubly-curved.


\end{abstract}

\section{Introduction}

Creases are lines of singularly large curvature formed when a thin sheet is folded beyond its elastic limit.  Regular patterns of creases follow controlled folding, where lines are typically straight; with care, they can be made to be curved, for example in thin paper card, by first scoring along predisposed directions.  Uncontrolled folding, such as occurs when a thin sheet is crumpled, results in a seemingly random layout of creases; but when the loads are removed and the sheet relaxes, the regions around the creases are developable, being locally flat or openly cylindrical or conical \cite{Pomeau}.  This also holds for controlled folding, and because simple curving is energetically efficient in engineering materials, it underpins the efficiency of, for example, sheet-metal working and the manufacturing of shaped panels and prismatic packaging.

Figure~1(a) shows a circular cylindrical tube made from a rectangular sheet of thin card, bent round and joined at the edges.  Before wrapping, a set of equally-spaced parallel lines is drawn onto the sheet at an angle $\alpha$ to the joining edges.  Defining a local coordinate system, $(x,y)$, with the $x$-axis parallel to the lines, we can easily compute the curvatures, $\kappa_{xx}$ and $\kappa_{yy}$, and twist, $\kappa_{xy}$, of the cylindrical surface with respect to $(x,y)$ as follows: $\kappa_{xx} = (1/a) \sin^2 \alpha$, $\kappa_{yy} = (1/a) \cos^2 \alpha$, and $\kappa_{xy} = (1/a) \sin \alpha \cos \alpha$, where $a$ is the cylinder radius.

Figure~1(b) shows what happens if, before wrapping, the lines on the sheet are scored with a sharp knife, and then these potential creases are all bent through a small angle, in order to make a uniform prismatic tube.  Now the edges of the sheet are not in register; and it is necessary to apply some shearing force to join them so that the edges of the tube are plane.  When that has been done, the tube looks superficially like Fig.~1(a); but now the strips between the lines no longer lie on the former cylindrical surface.  Instead, they are flat in the $y$-direction, so $\kappa_{yy}$ is zero; but otherwise, the values of $\kappa_{xx}$ and $\kappa_{xy}$ are the same as before.  In this ``twisted-prismatic'' form of the tube, the strips are no longer developable, as can be clearly seen in different views in Figs~1(c) and (d); they have non-zero Gaussian curvature; and the energy required to deform them by strain in the surface has been provided by the shearing force necessary to join up the edges.  An alternative way of achieving this final state would be to start with the circular cylindrical conformation and then, strip-by-strip, to push the material into the altered configuration with $\kappa_{yy} = 0$.  But that would be practically impossible to do, whilst the procedure described above achieves the same final result easily.  Importantly, the creased tube appears to contradict what we expect from controlled folding: that the strips between the creases are not simply curved.

In a previous paper \cite{KASCrease}, we showed by means of a Mohr's circle argument that the Gaussian curvature of the strips in the twisted-prismatic conformation is negative.  It is clear that a given length of tube subtends the same solid angle whether it is in the circular cylindrical or the twisted-prismatic form, and it follows that the total change of Gaussian curvature between the two states must be zero; hence, the increase in the Gaussian curvature in the curved creases separating the strips can be readily computed.  The purpose of the present note is to determine the Gaussian curvature of these creases by a more direct method: it is clear that any crease is characterised by its fold angle, the curvature of the crease line itself, and the twist along the crease, and we shall see how these three quantities affect and combine to express the required formula for the Gaussian curvature.

We first calculate the amount of negative Gaussian curvature in the strips using the details in \cite{KASCrease}.  We then focus attention directly on a crease and its solid angle performance via the related concepts of {\em spherical image} and {\em unit sphere}.  This is because the crease is infinitesimally thin,  of indeterminate area, so the Gaussian curvature is not properly defined.  We therefore investigate the performance of the solid angle by itself before defining a {\em specific} Gaussian curvature equal to the solid angle divided by the arc-length of crease.  This measure is independent of the width of surface in consideration around the crease: outside of the crease and in the strips, the same property is calculated by multiplying the Gaussian curvature of the strip by its width between successive creases.

We separate our study into two cases of a twisted crease and a curved crease.  We show that the solid angle of the former does not depend on the fold angle; we can unfold a twisted crease without changing the amount of Gaussian curvature.  On the other hand, the solid angle of the curved crease does depend on the fold angle, and we derive an expression for the specific Gaussian curvature in terms of this and the intrinsic curvature along it, which, when applied to the creased tube, verifies the balance we expect.  In two final exercises, we derive first the curved crease expression by modelling it as a continuous doubly-curved shell, and observe how the solid angle performs when the transverse width is reduced to zero.  In the second, we validate the expression separately for a creased shell whose surface is curved but not twisted.

\section{Analysis}

Each strip in the twisted-prismatic tube has non-zero, uniform $\kappa_{xx}$ and $\kappa_{xy}$, with $\kappa_{yy} = 0$.  The centre, $C$, of the corresponding Mohr's circle of curvature {\em vs} twist \cite{CRC} is located on the abscissa at $\kappa_{xx}/2 = (1/2a) \sin^2 \alpha$, and its radius, denoted as $B$, is equal to $\sqrt{\kappa_{xy}^2 + \kappa_{xx}^2/4}$, where

\begin{equation}
B^2 =  \biggl({\sin \alpha \cos \alpha \over a} \biggr)^2 + \biggl({\sin^2 \alpha \over 2 a} \biggr)^2
\end{equation}

\noindent The principal curvatures are given by $C \pm B$, and their product is $C^2 - B^2$, equal to the Gaussian curvature in the strip: multiplying this by the width normal to the creases, $h$, yields the specific Gaussian curvature in each strip as

\begin{equation}
h \left\{{\sin^4 \alpha \over 4a^2} - \left[{\sin^4 \alpha \over 4 a^2} + {\sin^2 \alpha \cos^2 \alpha \over a^2} \right]\right\} = -{h \over a^2} \sin^2 \alpha \cos^2 \alpha
\label{IGCB}
\end{equation}

\noindent This is clearly a negative quantity for all values of $\alpha$.

For each of the following crease scenarios, we calculate the solid angle by finding the area of the spherical image mapped onto a sphere of unit radius---the ``unit sphere'', as first proposed by Gauss \cite{CRC}.  Before giving the details of this mapping, the spherical coordinate system is outlined in Fig.~2.  The azimuthal angle around the equator is $\phi$ and the elevation angle relative to the equatorial plane is $\theta$.  A ring of surface width, $\delta \theta$, at a given elevation projects a horizontal radius equal to $\cos \theta$, and thus, has a surface area of $2 \pi \cos \theta \delta \theta$; for any portion of this ring subtending $\phi$ less than $2\pi$, the area is reduced by $\phi /2\pi$.

An element of a uniformly twisted surface is first shown in Fig.~3(a).  Even though we are determining {\em intrinsic} properties independent of any coordinates, we define an orthogonal system, $(x,y,z)$, with origin at the centre of the element for assistance in our calculations of solid angle.  The untwisted surface is the $(x,y)$ plane, and shallow twisted displacements normal to this are defined by $z = \kappa_{xy} x y$; the twist curvature, $\kappa_{xy}$, is therefore left-handed positive, which avoids the need for a minus sign in this expression.  The crease is eventually formed along the $x$-axis but, for now, we highlight a small rectangular patch inside with edges of length, $a$ and $b$, parallel to $x$ and $y$, respectively.

Unit vectors, ${\bf n}_1$, ${\bf n}_2$, ${\bf n}_3$ and ${\bf n}_4$, are drawn in sequence and normal to the surface at all corners of the patch.   These are then identically mapped onto the surface of the unit sphere in Fig.~3(b), in order to define the vertices of the corresponding spherical image.  This mapping preserves the orientation of the vector, which can be determined from the gradient field of the surface;  $\partial z / \partial x = \kappa_{xy} y$ {\em etc}.  However, along the edges of patch, the gradient only varies orthogonally to its direction, for example, along the $x$-edges, only $\partial z / \partial y$ changes. Therefore, the change in the orientation of the unit normals moving from corner to corner is $\kappa_{xy} a$ from ${\bf n}_1$ to ${\bf n}_2$, and $\kappa_{xy} b$ from ${\bf n}_2$ to ${\bf n}_3$, and so forth.  Furthermore, if we correlate the origin of $(x,y,z)$ with $(\phi, \theta) = (0,0)$ on the unit sphere, where $z$ is normal to the surface, the corners of the image are symmetrically located, and the edges form a spherical quadrilateral of side-lengths equal to $\kappa_{xy} a$ and $\kappa_{xy} b$, which we denote as $\delta \gamma$ and $\delta \xi$, respectively.

The cyclic ordering of the unit normals is also reversed on the unit sphere because of gradient orthogonality and, by convention, the area of the spherical image is negative.  This area, {\em viz.} the solid angle, is found from:

\begin{equation}
\int_{-\delta \gamma/2}^{-\delta \gamma/2} - 2\pi \cos \theta {\rm d} \theta \times {\delta \xi \over 2 \pi} = -\delta \xi \cdot 2 \sin{\delta \gamma \over 2}
\label{twistSA}
\end{equation}

\noindent which tends to $-\delta \xi \cdot \delta \gamma$ for shallow displacements and rotations; replacing each of these angles by $\kappa_{xy} a$ and $\kappa_{xy} b$, and dividing by the area of the original patch, $ab$, the familiar Gaussian curvature is recovered; $-\kappa_{xy}^2$.

A crease is now formed by folding the surface along the $x$-axis in Fig.~3(c) by equal and opposite fold angles, $\mu$.  Because the $x$-axis is straight, each half of the surface experiences a rigid body rotation without deformation, and there can be no change in the solid angle; creasing a shallow twisted surface is a developable process provided we crease along the principal axis of twist.  This is also verified when we consider how the spherical image of the twisted surface changes after creasing.

In Fig.~3(c), it is clear that each of the unit normals at the corners is tilted away from their original positions in Fig.~3(a).  Assuming that $\mu$ is also small, tilting takes place in the $y$-direction only and, on the unit sphere, their locations move by $\pm \mu$ along their original meridians.  Depending on which side they lie, all unit normals on the patch edges move by the same amount, and the outline of the original image is now distorted into an eight-sided polygon, see Fig.~3(d).  Its layout is antisymmetrical about the centre because there are step changes of $2 \mu$ in the orientation of unit normals moving across the crease at the front and back of the patch.  We may also think of the original spherical image shearing down the middle, as shown; and for small movements of each half towards either pole, it is clear that the total area remains the same.  This may be formally evaluated as follows:

\begin{equation}
\int_{-\delta \gamma/2+\mu}^{-\delta \gamma/2+\mu} - 2\pi \cos \theta {\rm d} \theta \times {(\delta \xi /2) \over 2 \pi} + \int_{-\delta \gamma/2-\mu}^{-\delta \gamma/2-\mu} - 2\pi \cos \theta {\rm d} \theta \times {(\delta \xi /2) \over 2 \pi} = -\delta \xi \cdot 2 \sin{\delta \gamma \over 2} \cos{\mu}
\label{twistCreaseSA}
\end{equation}

\noindent For small fold angles, $\cos \mu \approx 1$ and the result is the same as Eqn~\ref{twistSA}.  For moderate fold angles, where $\cos \mu < 1$, we cannot say there is a reduction in the solid angle because we violate the assumption of small rotations, as well as displacements, inherent in the construction of the spherical image.   We are satisfied that the current result stands, as does its corollary: that twisting along a crease does not change its solid angle.  We now consider the effect of curving along the crease.

The schematic view is given in Fig.~4, which shows the crease as a curved line of radius, $R$, connected on either side by thin strips of an arbitrary width.   The strips are inclined to the tangent plane of the crease by the same fold angle, $\mu$, and in this view, the arrangement resembles one axial element of the curved corrugated shells in \cite{KASCompliant}.  Simple handling of these shells show that they cannot be flattened without significant in-plane straining: in other words, they have inherent Gaussian curvature.  Although we are considering a static geometry without strains in Fig.~4, we expect the same to be true for a curved crease.

Figure~5(a) shows a patch symmetrically spanning the crease and subtending a small angle, $\delta \xi$ along it.   As before, normal unit vectors at the patch corners are highlighted for constructing the spherical image.  This is a relatively straightforward process if we again correlate $(\phi, \theta)$ to principal  directions respectively along the crease and normal to it within the side strips, with origin at the patch centre.  The image is, again, a spherical quadrilateral with side-lengths of $\delta \xi$ along $\phi$ and $2\mu$ along $\theta$: in Fig.~5(b) these have been exaggerated for clarity when the angles are small; and along $\theta$, we are capturing the abrupt change in the orientation of the unit normal moving across the crease.  The corner normals have the same cyclic order as on the patch, so the area is positive, which we denote  as $\delta \beta$, where

\begin{equation}
\delta \beta  = \int_{-\mu}^{+\mu} 2 \pi \cos \theta {\rm d} \theta \times {\delta \xi \over 2 \pi} = \delta \phi \cdot 2 \sin \mu
\label{dBdyDef}
\end{equation}

\noindent This clearly depends on the fold angle but it is independent of the width of patch: no matter how narrow the patch becomes, the solid angle remains finite and thus, we conclude that the Gaussian curvature of the crease lies entirely {\em within} the crease itself.  Such singular behaviour is not uncommon \cite{CRC}, for example, in discrete roof vertices and, more generally, in polyhedral approximations of generally curved surfaces where often, some definition of associated area is often brought to bear, in order to precisely compute the overall shape of structure.  Here, we dispense with area and divide Eqn~\ref{dBdyDef} by the arc-length of crease, akin to the specific Gaussian curvature defined previously.  Denoting $\delta x$ as $R \delta \phi$, where $x$ is the usual coordinate along the crease, we find $\delta \beta = 2 (\delta x /R) \sin \mu$.  Dividing by $\delta x$ and observing the limit, we obtain:

\begin{equation}
{{\rm d} \beta \over {\rm d} x} = {2 \sin \mu \over R}
\label{dBdy}
\end{equation}

\noindent which may be expressed more generally as:

\begin{equation}
{2 \sin \mu \over R} = 2 \times \text{crease curvature along} \times \sin{(\text{0.5 $\times$ total fold angle across})}
\label{IGC}
\end{equation}

\noindent Now consider the creased tube of Fig.~1, where the ``curvature along'' is simply $\kappa_{xx} = (1/a) \sin^2 \alpha$.  The total fold angle is equivalent to the angle subtended between the lines drawn onto the smooth cylinder in Fig.~1(a), itself equal to $\kappa_{yy}$ multiplied by the normal width, $h$.  Substituting these values into the above returns

\begin{equation}
2 \times {\sin^2 \alpha \over a} \times \sin \biggl({h \over 2a} \cos^2 \alpha \biggr) \approx {h \over a^2} \sin^2 \alpha \cos^2 \alpha
\label{smallIGC}
\end{equation}

\noindent after assuming that $\sin \mu \approx \mu$ for small $\mu$.  Compared to Eqn~\ref{IGCB}, this expression is equal and opposite, and the balance of Gaussian curvature throughout the tube is confirmed.

We now consider two final examples of separate emphasis.  The first replicates Eqn~\ref{dBdy} using an independent derivation based on an {\em extrinsic} viewpoint---by treating the crease initially as a continuous surface: we do not consider, nor need, a spherical image, rather we operate directly with the surface specification.  In the second, we apply Eqn~\ref{dBdy} to a different shell made of pieces connected at seams that perform as curved creases; there is no folding of the surface during construction but Eqn~\ref{dBdy} is equally valid.

First, the crease is modelled as a narrow, doubly-curved surface similar to a bicycle mudguard, as shown in Fig.~6(a).  The surface exactly inscribes the original discrete set-up from Fig.~4, where the edges of both cross-sections touch tangentially.  The surface is formed when a transverse arc, of radius $r$, is rotated around an axis at the much larger radius, $R$, of the discrete crease.  One principal radius of curvature is $r$ and the other, in the swept direction, is $[R - r(1- \cos \mu)]/\cos \epsilon$, which applies to any point on the inscribed arc at the local angle, $\epsilon$, in Fig.~4(b).  The principal curvatures are simply reciprocal expressions, and their product defines the local Gaussian curvature, which may be integrated over a complete hoop to yield the total solid angle of the surface:

\begin{equation}
2 \pi R \times \int_{-\mu}^{+\mu} {1 \over r} \cdot {\cos \beta \over R - r(1 - \cos \mu)} r {\rm d} \beta = {4 \pi R \sin \mu \over R - r (1 - \cos \mu)}
\end{equation}

\noindent As $r$ vanishes, the surface approaches a discrete crease, obviating $R >> r$ in the above and leading to $4 \pi \sin \mu$.  Dividing by $2 \pi R$, we obtain the ratio of solid angle to arc-length, and the final expression is the same as Eqn~\ref{dBdy}.

Second, consider the approximation to a sphere in Fig.~7.  This is formed by bending originally flat, thin tapering pieces of plate along their meridians and then joining them along adjacent seams.  The precise shape of pieces, or gores, is key because this enables construction of the sphere to be developable without in-plane distortion or tearing of the gores; we shall assume this to be the case. But it is clear that the sphere has Gaussian curvature in a average sense, which tends to $1/R^2$ when the number of gores, $n$, becomes large enough, where $R$ is the radius of the seam; and in this case, the solid angle subtended by the discrete shape tends to $4\pi$, the standard result for a smooth sphere.

Each of the gores is cylindrical, and the Gaussian curvature is entirely contained within the seams, which we assume to behave as curved creases.  At a given elevation, the latitudinal width, $w$, between each pair of creases is $2\pi r/n$, where $r = R \cos \theta$ as shown.  The curvature of the latitudinal circle is $1/r_2$ by definition, equal to $\cos \theta /R$; when this is multiplied by $w$, we obtain the total fold angle normal to the crease, that is $2 \mu = w \cos \theta /r$.  Upon replacing $w$ and $r$, then $\mu = (\pi / n) \cos \theta$, and setting the elemental arc-length, ${\rm d} x$, equal to $R {\rm d} \theta$, Eqn~\ref{dBdy} produces

\begin{equation*}
{{\rm d} \beta \over {\rm d} \theta} = 2 \sin\left({\pi \over n} \cos \theta \right)
\end{equation*}

\noindent for a single crease. For a large number of creases, the fold angle is everywhere small, which enables the above to be rewritten as

\begin{equation*}
{{\rm d} \beta \over {\rm d} \theta} \approx 2 {\pi \over n} \cos \theta~~\rightarrow~~ n \times \int_{-\pi/2}^{\pi/2} {\rm d} \beta = 2 \pi \int_{-\pi/2}^{\pi/2} \cos \theta {\rm d} \theta = 4 \pi
\end{equation*}

\noindent as per the smooth spherical case.




\section{Conclusions}

In an earlier paper we observed that the strips between creases in a creased tube have negative Gaussian curvature but we did not study this property any further: we did not realise that there must be a zero net result and that the balance is created by the creases themselves.  Here, we have computed the Gaussian curvature of the crease in terms of the performance of its solid angle with respect to the crease arc-length. This approach is necessary because the crease is a discrete feature without width, where the conventional definitions of Gaussian curvature are not meaningful. Our calculation of  specific Gaussian curvature in the creases exactly balances that of the strips, provided the fold angle across the crease is small.  The performance also depends on the curvature along the crease but not on any surface twist.  The latter point may be surprising to some because it is often assumed that a surface with inherent Gaussian curvature can only deform extensionally: this is always true when there is positive Gaussian curvature; but as we have seen, not always so when it is negative.  Our analysis has made use of the spherical image because of the relatively simple, uniform surface geometry around the crease when the fold angle is small, and we have verified our relationship from a continuum viewpoint.

~

\newpage
\section*{Figures}

\psfrag{m}{$\mu$}
\psfrag{2m}{$2\mu$}
\psfrag{B}{\hspace{-0.25cm} $\beta$}
\psfrag{n1}{\bf{n}$_1$}
\psfrag{n2}{\bf{n}$_2$}
\psfrag{n3}{\bf{n}$_3$}
\psfrag{n4}{\bf{n}$_4$}
\psfrag{n1}{\bf{n}$_1$}
\psfrag{R}{$R$}
\psfrag{r}{$r$}
\psfrag{x}{$x$}
\psfrag{y}{$y$}
\psfrag{z}{$z$}
\psfrag{p}{$\rho$}
\psfrag{axis}{axis}
\psfrag{df}{$\delta \phi$}
\psfrag{(a)}{(a)} \psfrag{(b)}{(b)} \psfrag{(c)}{(c)}
\psfrag{(d)}{(d)} \psfrag{(e)}{(e)}

\begin{figure}[!ht]
\begin{center}
\includegraphics[width=16cm]{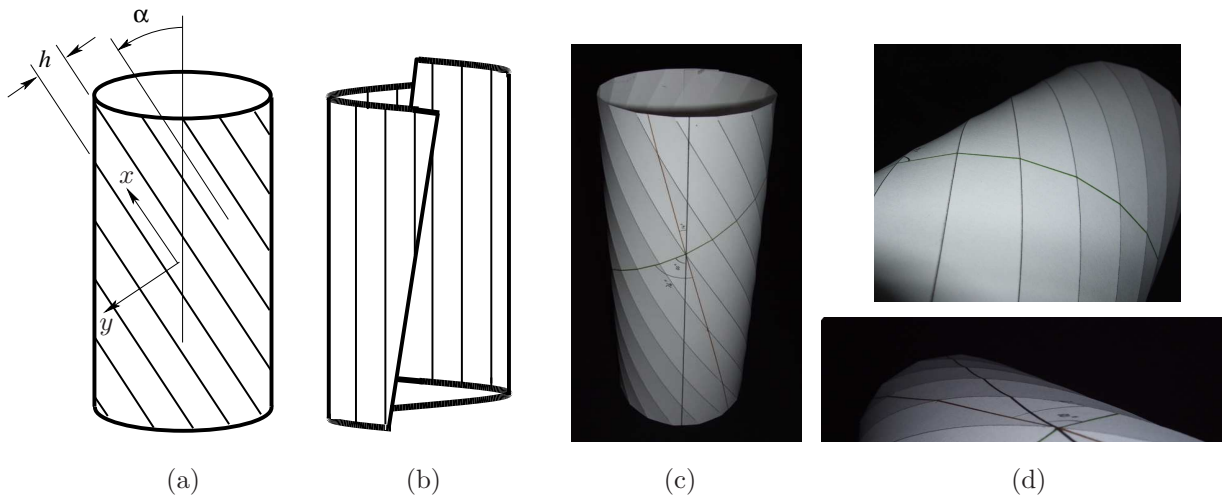}
\caption{Manufacture and shape of a creased tube from \cite{KASCrease}. (a) Parallel lines are drawn onto flat sheet of card inclined at angle, $\alpha$, to one edge, before wrapping into a cylinder.  A local orthogonal coordinate system, $(x,y)$, is embedded in the surface with $x$ parallel to the lines.  (b) The lines in (a) are now scored with a knife before wrapping, giving sharp, folded creases.  Overall, the stress-free shape is helicoidal, where the strips between creases remain flat.  (c) The free edges in (b) are now joined together to form a twisted-prismatic tube made of thin card.  The strips remain flat normal to creases, but there is a second direction of flatness highlighted in this example, whose inclination was calculated in \cite{KASCrease}. (d) Close-up views showing a flat line normal to creases, top, and a gentle, inwardly curved axial line, signifying the double curvature, bottom.}
\end{center}
\end{figure}

\begin{figure}[!ht]
\begin{center}
\psfrag{p}{$\phi$} \psfrag{q}{$\theta$} \psfrag{dq}{$\delta \theta$} \psfrag{O}{$O$} \psfrag{P}{$P$}
\includegraphics[width=6cm]{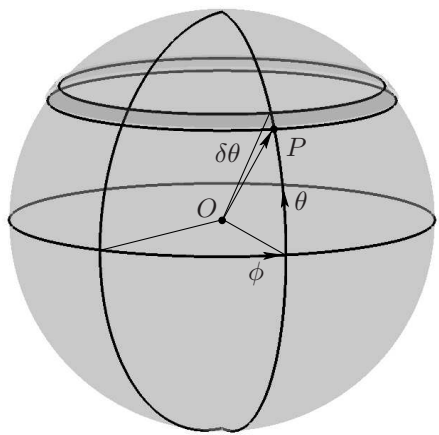}
\caption{Spherical coordinate system for calculating the area of a patch or ring on the surface of a unit sphere: $\phi$ is the azimuthal angle and $\theta$ is the elevation angle above the equator.  The origin, $O$, lies at the centre of sphere, and a point, $P$ is highlighted on a ring of width $\delta \theta$.  The surface area of this ring is $2 \pi \cos{\theta} \delta \theta$.}
\end{center}
\end{figure}

\begin{figure}[!ht]
\begin{center}
\psfrag{2b}{$\delta \xi$} \psfrag{2a}{$\delta \gamma$}
\includegraphics[width=11cm]{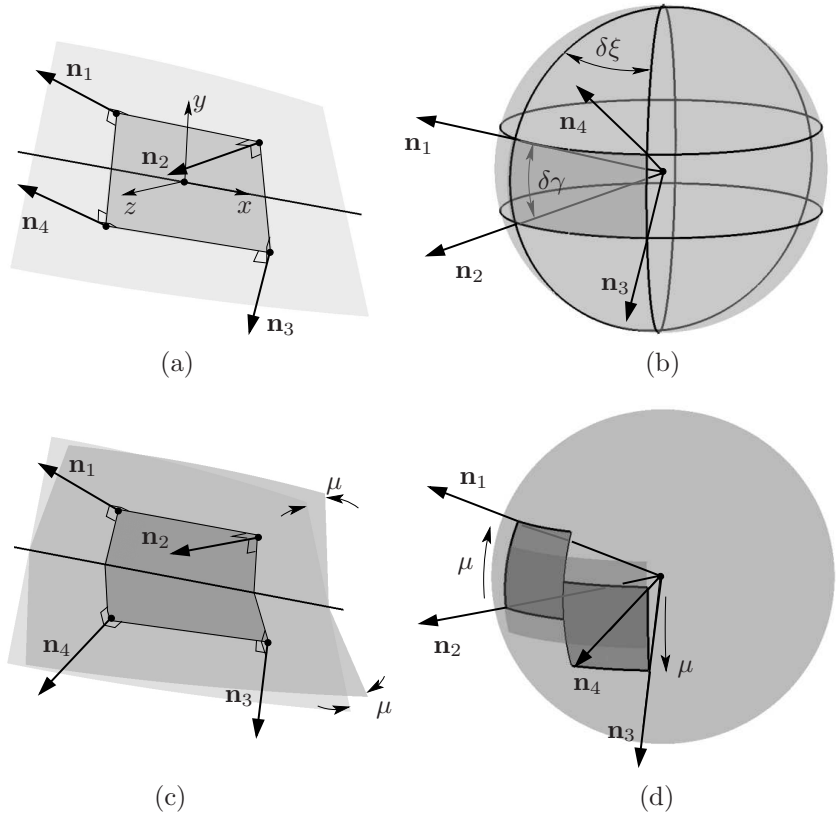}
\caption{Calculation of the solid angle of a twisted crease. (a) Ordinary twisted surface with local coordinate system, $(x,y,z)$, located at the centre of an inner rectangular patch of side-lengths $a$ and $b$ along $x$ and $y$: the surface is described by $z = \kappa_{xy} x y$ relative to the $(x,y)$ plane, with the twisting curvature equal to $\kappa_{xy}$, and unit vectors, ${\bf n}_1 \ldots {\bf n}_4$ are drawn normal to the surface of the patch at its corners. (b) Spherical image of the patch in (a) after mapping the same normals identically onto the surface of the unit sphere.  The image is quadrilateral and subtends $\delta \xi = \kappa_{xy} b$ and $\delta \gamma = \kappa_{xy} a$ in the $\phi$ and $\theta$ directions, respectively.  The image area equates to the magnitude of the solid angle of the patch, and the cyclic direction of normals is opposite to that in (a), signifying a negative area.  (c) Twisted crease where the surface in (a) has been folded about the $x$-axis by $\pm \mu$.  (d) Spherical image of the patch in (c) after conserving the orientation of unit normals: the original image is shown behind, in order to convey the effect of creasing.}
\end{center}
\end{figure}

\begin{figure}[!ht]
\begin{center}
\includegraphics[width=8cm]{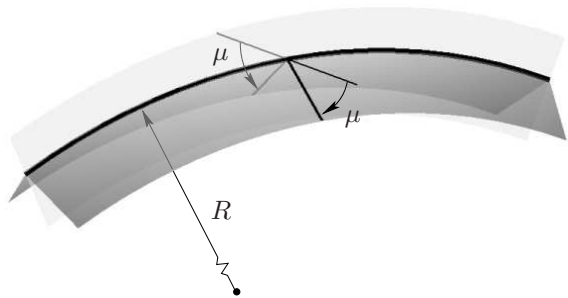}
\caption{A curved crease  modelled as an arc (black) of radius, $R$, connected to conical strips.  These are inclined at equal and opposite fold angles, $\mu$, everywhere relative to the tangent plane of the crease.  The plane of $R$ bisects the total fold angle.}
\end{center}
\end{figure}

\begin{figure}[!ht]
\begin{center}
\psfrag{df}{$\delta \xi$}
\includegraphics[width=14cm]{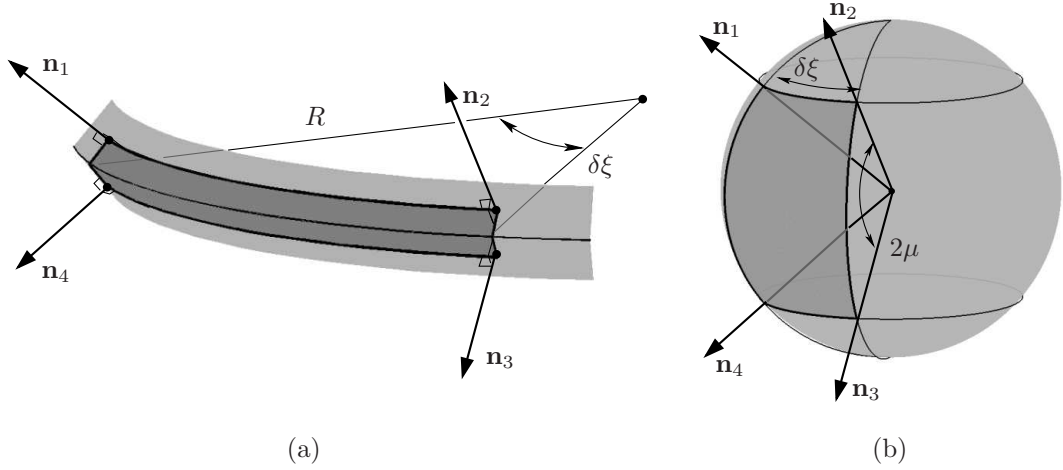}
\caption{Calculation of the solid angle of a curved crease.  (a) Unit normals are highlighted at the corners of a rectangular patch encapsulating part of the crease.  (b) Spherical image of the patch in (a), where the edges subtend $\delta \xi$ and $2\mu$.  The ordering of normals is the same as in (a) and the image area is positive.}
\end{center}
\end{figure}

\begin{figure}[!ht]
\begin{center}
\psfrag{B}{\hspace{-0.25cm} $\epsilon$}
\includegraphics[width=12cm]{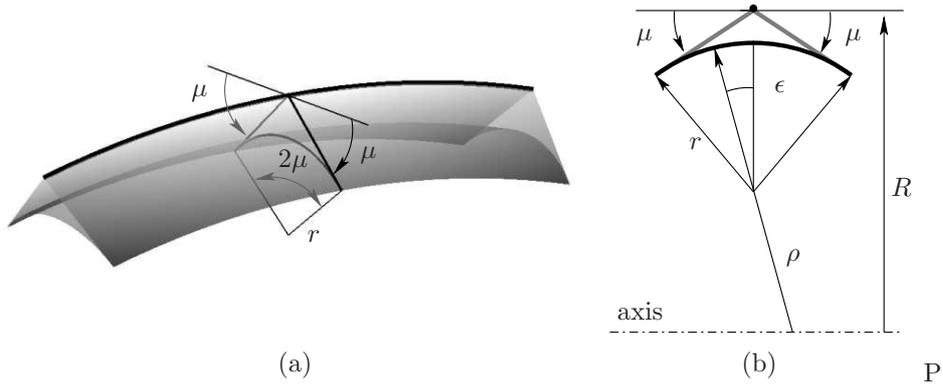}
P\caption{Continuous approximation of a curved crease.  (a) The crease is rendered as a doubly-curved surface exactly inscribing the discrete version from Fig.~4. The cross-section subtends $2\mu$, the total fold angle, at radius, $r$. (b) The surface in (a) is formed by sweeping the transverse arc around an axis at radius, $R$, of the original discrete crease. $\rho$ is the radius of curvature of the surface in the swept direction, and $\epsilon$ is an internal coordinate on the cross-section.}
\end{center}
\end{figure}

\begin{figure}[!ht]
\begin{center}
\psfrag{r2}{$r_2$} \psfrag{w0}{$w_0$} \psfrag{ph}{} \psfrag{p/n}{$\pi/n$} \psfrag{2p/n}{$2 \pi / n$} \psfrag{w}{$w$} \psfrag{q}{$\theta$}
\includegraphics[width=6cm]{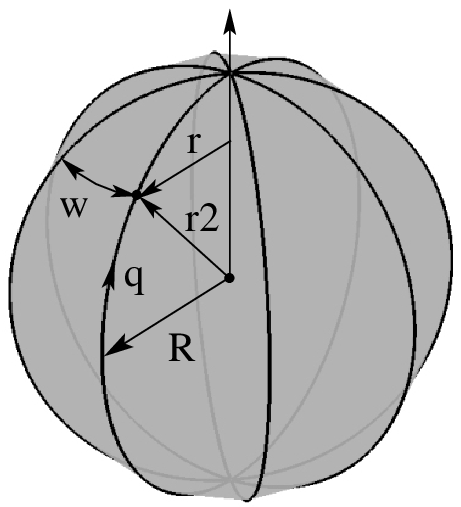}
\caption{Discrete approximation of a sphere in which $n$ tapering gores are connected on their edges to form seams with radius of curvature, $R$.  The latitudinal width at a given elevation is $w$ and the associated radius of curvature is $1/r_2$ equal to $1/R$ for a sphere.}\end{center}
\end{figure}


\begin{thebibliography}{}

\bibitem{Pomeau}
M Ben Amar and Y Pomeau, Crumpled paper, {\em Proceedings of the Royal Society of London, Series A},  453 (1997), pp.~729-755

\bibitem{KASCrease}
K A Seffen, N Borner, The shape of helically creased cylinders, {\em Journal of Applied Mechanics, ASME} (2013), 80(5)

\bibitem{CRC}
C R Calladine, {\em Theory of Shell Structures} (1983), Cambridge University Press



\bibitem{KASCompliant}
K A Seffen, Compliant shell mechanisms, {\em Philosophical Transactions of the Royal Society of London A}, (2012) 370 pp.~2010-2026







\end{thebibliography}
\end{document}